\documentclass[12pt,showkeys]{article}
\usepackage{amssymb}

\title{Occupation densities for SPDE's with
Reflection}
\author{\Large{Lorenzo Zambotti}
\\ Scuola Normale Superiore\\
Piazza dei Cavalieri 7, 56126 Pisa, Italy
\\ E-mail: zambotti@sns.it} 
\date{}

\newtheorem{proposition}{Proposition}[section]
\newtheorem{theorem}{Theorem}[section]
\newtheorem{lemma}{Lemma}[section]
\newtheorem{corollary}{Corollary}[section]
\newtheorem{definition}{Definition}[section]
\newtheorem{remark}{Remark}[section]

\begin{document}

\maketitle

\begin{abstract}
We consider the solution $(u,\eta)$ of the
white-noise driven stochastic partial differential
equation with reflection on the space interval $[0,1]$
introduced by Nualart and Pardoux.
First, we prove that at any fixed time $t>0$, the measure
$\eta([0,t]\times d\theta)$ is absolutely continuous w.r.t.
the Lebesgue measure $d\theta$ on $(0,1)$. We characterize the
density as a family of additive functionals of $u$, and we 
interpret it as a renormalized local time
at $0$ of $(u(t,\theta))_{t\geq 0}$. Finally we study the
behaviour of $\eta$ at the boundary of $[0,1]$. The main technical
novelty is a projection principle from the Dirichlet space
of a Gaussian process, vector-valued solution of a linear
SPDE, to the Dirichlet space of the process $u$.
\end{abstract}

\section{Introduction}

We are concerned with the solution $(u,\eta)$ of the
stochastic partial
differential equation with reflection 
of the Nualart-Pardoux type, see \cite{nupa}:
\begin{equation}\label{1}
\left\{ \begin{array}{ll}
{\displaystyle
\frac{\partial u}{\partial t}=\frac 12
\frac{\partial^2 u}{\partial \theta^2}
 + \frac{\partial^2 W}{\partial t\partial\theta} +
\eta(t,\theta) }
\\ \\
u(0,\theta)=x(\theta), \ u(t,0)=u(t,1)=0
\\ \\
u\geq 0, \ d\eta\geq 0, \
\int u\, d\eta =0
\end{array} \right.
\end{equation}
where $u$ is a continuous
function of $(t,\theta)\in {\overline{\cal O}}:=
[0,+\infty)\times[0,1]$,
$\eta$ a positive measure on ${\cal O}:=
[0,+\infty)\times(0,1)$, $x:[0,1]\mapsto[0,\infty)$
and $\{W(t,\theta):\ (t,\theta)\in
{\overline{\cal O}}\}$ is a Brownian sheet. We denote
by $\nu$ the law of a Bessel Bridge 
$(e_\theta)_{\theta\in[0,1]}$ of dimension 3.

The main aim of this paper is to prove the following properties
of the reflecting measure $\eta$:
\begin{itemize}
\item[1.] For all $t\geq 0$ the measure $\eta([0,t],d\theta)$ 
is absolutely continuous with respect to the Lebesgue measure
$d\theta$ on $(0,1)$:
\begin{equation}\label{abscon}
\eta([0,t],d\theta) \, = \, \eta([0,t],\theta) \, d\theta.
\end{equation}
The process $(\eta([0,t],\theta))_{t\geq 0}$, $\theta\in(0,1)$, 
is an Additive Functional of $u$, increasing only 
on $\{t:u(t,\theta)=0\}$, with Revuz measure:
\begin{equation}\label{rev1}
\frac 1{2\sqrt{2\pi\theta^3
(1-\theta)^3}} \, \nu( \, dx \, | \, x(\theta)=0).
\end{equation}
\item[2.] For all $t\geq 0$:
\begin{equation}\label{intl1}
\eta([0,t],\theta) \, = \, \frac 34
\, \lim_{\epsilon\downarrow 0}
\frac 1{\,\epsilon^3} \int_0^t 1_{[0,\epsilon]}(u(s,\theta)) \, ds,
\end{equation}
in probability.
\item[3.] There exists a family of Additive Functionals of
$u$, $(l^a(\cdot,\theta))_{a\in[0,\infty), \,
\theta\in(0,1)}$, such that $l^a(\cdot,\theta)$ increases
only on $\{t:u(t,\theta)=a\}$ and such that the following
occupation times formula holds for all
$F\in B_b({\mathbb R})$: 
\begin{equation}\label{occui}
\int_0^t F(u(s,\theta)) \, ds \, = \, 
\int_0^\infty F(a) \, l^a(t,\theta) 
\, da, \quad t \geq 0.
\end{equation}
\item [4.] For all $t\geq 0$:
\begin{equation}\label{intl4}
\eta([0,t],\theta) \, = \, \frac 14 \,
\lim_{a\downarrow 0} \frac 1{\, a^2} \, l^a(t,\theta)
\end{equation}
in probability.
\item[5.] For all $t\geq 0$ and $a\in(0,1)$:
\begin{equation}\label{intl3}
\lim_{\epsilon\downarrow 0}\, \sqrt\epsilon \,
\int_0^a \left(1 \wedge \frac{\theta}\epsilon \right)
\, \eta([0,t], d\theta) \, = \, {\sqrt{\frac 2\pi}} \, t 
\end{equation}
and symmetrically:
\begin{equation}\label{intl3'}
\lim_{\epsilon\downarrow 0}\, \sqrt\epsilon \,
\int_a^1 \left(1 \wedge \frac{1-\theta}\epsilon \right)
\, \eta([0,t], d\theta) \, = \, {\sqrt{\frac 2\pi}} \, t,
\end{equation}
in probability.
\end{itemize}
Recall that if $B$ is a linear Brownian Motion and $(X,L)$
is the unique continuous solution of the Skorohod problem:
\[
dX \, = \, dB \, + \, dL, \quad X(0)=x\geq 0, \ L(0)=0,
\]
\[
X\geq 0, \quad t\mapsto L(t) \ {\rm non-decreasing},
\quad \int_0^\infty X(t) \, dL(t) \, = \, 0,
\]
then it turns out that $2L$ is the local time of $X$ at
$0$ and:
\begin{equation}\label{refbm}
L(t) \, = \, \frac 12 \, \lim_{\epsilon\to 0} \, \frac
1\epsilon \int_0^t 1_{[0,\epsilon]}(X(s)) \, ds.
\end{equation}
In the infinite-dimensional equation (\ref{1}), the
reflecting term $\eta$ is a random measure on space-time. 
In \cite{za00}, the following decomposition formula was 
proved:
\begin{equation}\label{decom00}
\eta(ds,d\theta) \, = \, \delta_{r(s)}(d\theta) \, 
\eta(ds,(0,1)),
\end{equation}
where $\delta_a$ is the Dirac mass at $a\in(0,1)$ and 
$r(s)\in(0,1)$, for $\eta(ds,(0,1))$-a.e. $s$, is
the unique $r\in(0,1)$ such that $u(s,r)=0$. This formula
was used in \cite{za00} to write equation (\ref{1}) as
the following Skorohod problem in the infinite dimensional convex set 
$K_0$ of continuous non-negative $x:[0,1]\mapsto[0,\infty)$:
\[
du \, = \, \frac 12 \frac {\partial^2 u}{\partial \theta^2}
 \, dt \, + \, dW \, + \, \frac 12 \, n(u) \cdot dL,
\]
interpreting the set of $x\in K_0$ having a unique zero
in $(0,1)$ as the boundary of $K_0$, the increasing 
process $t\mapsto L_t:=2\eta([0,t],(0,1))$ as the local time
of $u$ at this boundary and the measure
$n(u)=\delta_{r(s)}$ as the normal vector field to this boundary
at $u(s,\cdot)$.

On the other hand,
the absolute-continuity result (\ref{abscon}) suggests
an interpretation of $\eta$ as sum of
reflecting processes $t\mapsto \eta([0,t],\theta)$, each
depending only on $(u(t,\theta))_{t\geq 0}$ and 
increasing only on $\{t:u(t,\theta)=0\}$. Therefore, by
(\ref{abscon}) equation (\ref{1}) can also be interpreted 
as the following infinite system of 1-dimensional
Skorohod problems, parametrized by
$\theta\in(0,1)$ and coupled through the interaction given
by the second derivative w.r.t. $\theta$:
\begin{equation}\label{sk11}
\left\{ \begin{array}{ll}
{\displaystyle
u(t,\theta)=x(\theta) + \frac 12 \int_0^t 
\frac{\partial^2 u}{\partial \theta^2}(s,\theta) \, ds 
+ \frac{\partial W}{\partial \theta}(t,\theta)
+ \eta([0,t],\theta)
}
\\ \\
u(t,0)=u(t,1)=0
\\ \\
{\displaystyle
u\geq 0, \ \eta(dt,\theta)\geq 0, \
\int_0^\infty u(t,\theta)\, \eta(dt,\theta) =0, \quad \forall
\theta\in(0,1), }
\end{array} \right.
\end{equation}
see (\ref{11}) below.
This interpretation is reminiscent of the result of Funaki and
Olla in \cite{fuol}, where the fluctuations around the
hydrodynamic limit of a particle system with reflection
on a wall is proved to be governed by the SPDE (\ref{1}).

By (\ref{occui}), $(u(t,\theta))_{t\geq 0}$ admits for all $a\geq 0$
a local time at $a$, $(l^a(t,\theta))_{t\geq 0}$. However, 
by (\ref{intl4}), the reflecting term
$\eta([0,\cdot],\theta)$ which appears in (\ref{sk11}) is
not proportional to $l^0(\cdot,\theta)$, which in fact turns out
to be identically $0$, and is rather a renormalized local time.
The necessity of such renormalization is linked with the
unusual rescaling of (\ref{intl1}). These two properties
of $\eta$ seem to be significant differences 
w.r.t. the finite-dimensional Skorohod problems.

The formulae (\ref{intl3}) and (\ref{intl3'}) give informations
about the behaviour of $\eta$ near the boundary of $[0,1]$. 
In particular, (\ref{intl3}) and (\ref{intl3'}) prove that for
any $t>0$ and any initial condition $x$, the mass of $\eta$
on $[0,t]\times(0,1)$ is infinite. This solves a problem
posed by Nualart and Pardoux in \cite{nupa}. Notice also 
that the right hand sides of (\ref{intl3})-(\ref{intl3'}) 
are independent of the initial condition $x$.

\vspace{.3cm}
In \cite{za00} it was proved that for all
$I\subset\!\subset (0,1)$, the process 
$t \, \mapsto \, \eta([0,t]\times I)$, where $\eta$ is
the reflecting term of (\ref{1}), is an Additive Functional
of $u$, with Revuz measure:
\begin{equation}\label{revin'}
\frac 12 \, 
\int_I \frac 1{\sqrt{2\pi\theta^3(1-\theta)^3}} \, 
\, \nu(\, dx \,| \, x(\theta)=0) \, d\theta.
\end{equation}

At a heuristic level, the informations given by the 
formulae (\ref{abscon}), 
(\ref{intl1}), (\ref{intl4}), (\ref{intl3}) and (\ref{intl3'})
are already contained in (\ref{revin'}) and in the properties 
of the invariant measure $\nu$ of $u$: for instance, if
the limit in the right-hand side of (\ref{intl1}) exists
for all $\theta\in(0,1)$, then by the properties of $\nu$
the Revuz-measure of the limit is (\ref{rev1}) and therefore
(\ref{abscon}) holds by (\ref{revin'}) and by the 
injectivity of the Revuz-correspondence.

However, the existence of such limit 
is not implied by the structure of (\ref{revin'}) alone.
According to the Theory of Dirichlet Forms, a sufficient
condition for the convergence
of a family of additive functionals of a Markov process, as
for instance in (\ref{intl1}), is the convergence
in the Dirichlet space of the corresponding 1-potentials:
see Chapter 5 of \cite{fot}. In
our case, this amounts to introduce the potentials:
\[
U_\epsilon(x) \, := \, \frac 34
\, \int_0^\infty e^{-t} \,
\frac 1{\,\epsilon^3} \, {\mathbb E}\left[ 
1_{[0,\epsilon]}(u(s,\theta))\right] \, ds,
\]
where $x:[0,1]\mapsto[0,\infty)$ is continuous and $u$ is
the corresponding solution of (\ref{1}), and prove that
$U_\epsilon$ has a limit as $\epsilon\to 0$ with
respect to the Dirichlet Form:
\[
{\cal E}(\varphi,\psi) \, := \,
\frac 12 \, \int \langle \nabla \varphi,
\nabla \psi \rangle \, d\nu, \qquad \varphi,\psi\in
W^{1,2}(\nu),
\]
where $\nabla$ and $\langle\cdot,\cdot\rangle$ denote
respectively the gradient and the canonical scalar product
in $H:=L^2(0,1)$. 
Indeed, as proved in \cite{za00}, $u$ is the diffusion
properly associated with ${\cal E}$ in $L^2(\nu)$.

However, due to the strong irregularity of the reflecting
measure $\eta$ in (\ref{1}), a direct computation of the
norm of the gradient of $U_\epsilon$ seems to be out of
reach. In order to overcome this difficulty, we take
advantage of a connection between equation (\ref{1})
and the following ${\mathbb R}^3$-valued linear SPDE
with additive white-noise:
\begin{equation}\label{ou3}
\left\{ \begin{array}{ll}
{\displaystyle
\frac{\partial z_3}{\partial t}=\frac 12
\frac{\partial^2 z_3}{\partial \theta^2}
 + \frac{\partial^2 \overline{W}_3}{\partial t\partial\theta} }
\\ \\
z_3(t,0)=z_3(t,1)=0
\\ \\
z_3(0,\theta)\, =\overline{x}(\theta) 
\end{array} \right.
\end{equation}
where $\overline{x}\in H^3$ and 
$\overline{W}_3$ the is the ${\mathbb R}^3$-valued
Gaussian process whose components
are $3$ independent copies of $W$.
The process
$z_3$ is also called the ${\mathbb R}^3$-valued random
string (see \cite{fun} and \cite{mutr}), 
and is the diffusion properly associated with the 
Dirichlet Form in $L^2(\mu_3)$:
\[
\Lambda^3(F,G):=\int_{H^3} \langle \overline{\nabla} F
, \overline{\nabla} G \rangle_{H^3} \,
d\mu_3
\quad F,G\in W^{1,2}(\mu_3),
\]
where $\mu_3$ is the law in $H^3$ of a standard ${\mathbb
R}^3$-valued Brownian Bridge, and
$\overline{\nabla} F:H^3\mapsto H^3$ is the gradient
of $F$ in $H^3$. Then, in \cite{za00} it was noticed that
the Dirichlet Form ${\cal E}$ is the image of $\Lambda^3$
under the map $\Phi_3:H^3\mapsto H$, 
\ $\Phi_3(y)(\theta):=|y(\theta)|_{{\mathbb R}^3}$, i.e. 
$\nu$ is the image $\mu_3$ under $\Phi_3$ and:
\[
W^{1,2}(\nu)=\{\varphi\in L^2(\nu)
: \varphi\circ\Phi_3\in W^{1,2}(\mu_3)\},
\]
\[
{\cal E}(\varphi,\psi) \, = \, \Lambda^3
(\varphi\circ\Phi_3,\psi\circ\Phi_3) \quad
\forall \, \varphi,\psi \in W^{1,2}(\nu).
\]
This connection involves directly the Dirichlet Forms ${\cal
E}$ and $\Lambda^3$, but not the corresponding processes.
In particular, it does not imply that $u$ is equal in law to
$|z_3|$. Nevertheless, in this paper we prove that this
connection gives a useful projection principle from 
$W^{1,2}(\mu_3)$ onto $W^{1,2}(\nu)$ and that, in particular,
the convergence in $W^{1,2}(\mu_3)$ of the 1-potentials of $z_3$: 
\[
\overline{U}_\epsilon(\overline{x}) \, := \,
\, \frac 34 \, \int_0^\infty e^{-t} \,
\frac 1{\,\epsilon^3} \, {\mathbb E}\left[ 
1_{[0,\epsilon]}(|z_3(s,\theta)|)\right] \, ds,
\]
as $\epsilon\to 0$, implies the convergence of
the 1-potentials $U_\epsilon$ of $u$ in
$W^{1,2}(\nu)$, and therefore that (\ref{intl1}) holds.
Also the formulae (\ref{intl4}), (\ref{intl3}) and
(\ref{intl3'}) are proved similarly. Therefore, precise and
non-trivial informations about $u$ can be obtained from
the study of the Gaussian process $z_3$. 

We recall that an analogous
connection has been proved in \cite{zaza} to hold 
between the ${\mathbb R}^d$-valued
solution of a linear white-noise driven SPDE, $d\geq 4$, and the
solution of a real-valued non-linear white-noise driven 
SPDE with a singular drift.

\vspace{.3cm}
The paper is organized as follows. Section 2 contains the
main definitions and the preliminary results on potentials of
the random string in dimension 3. In section 3 the occupation 
densities and the occupation times formula (\ref{occui})
are obtained for the SPDE with reflection (\ref{1}). 
The main results, together with some
corollaries, are then proved in section 4.

\section{The 3-dimensional random string}

We denote by $(g_{t}(\theta,\theta'):t>0, \
\theta,\theta'\in(0,1))$ the fundamental solution of the
heat equation with homogeneous Dirichlet boundary condition,
i.e.:
\[
\left\{ \begin{array}{ll}
{\displaystyle
\frac{\partial g}{\partial t}=\frac 12
\frac{\partial^2 g}{\partial \theta^2}
}
\\ \\
g_t(0,\theta')=g_t(1,\theta')=0
\\ \\
g_0(\theta,\cdot)=\delta_\theta
\end{array} \right.
\]
where $\delta_a$ is the Dirac mass at $a\in(0,1)$.
Moreover, we set $H:=L^2(0,1)$ with the canonical scalar
product $\langle \cdot,\cdot \rangle$ and norm $\|\cdot\|$, 
$K_0:=\{x\in H:x\geq 0\}$, 
\[
C_0:=C_0(0,1):= \{c:[0,1]\mapsto{\mathbb R}\ \
{\rm continuous},\ c(0)=c(1)=0\},
\]
\[
A:D(A)\subset H\mapsto H, \quad D(A):=W^{2,2}\cap W^{1,2}_0(0,1), \quad
A:=\frac 12 \frac{d^2 }{d \theta^2},
\]
and for all Fr\'echet differentiable $F:H\mapsto{\mathbb R}$
we denote by $\nabla F:H\mapsto H$ the gradient in $H$.
We set ${\cal O}:=[0,+\infty)\times(0,1)$ and 
${\overline{\cal O}}:=[0,+\infty)\times[0,1]$.
We denote by $(e^{tA})_{t\geq 0}$ the
semigroup generated by $A$ in $H$, i.e.:
\[
e^{tA}h(\theta) \, := \, \int_0^1 g_t(\theta,\theta') \,
h(\theta') \, d\theta', \qquad h\in H.
\]
Let $W$ be a two-parameter 
Wiener process defined on a complete probability space 
$(\Omega,{\cal F}, {\mathbb P})$, i.e. a 
Gaussian process with zero mean and covariance function
\[
{\mathbb E}\, \left[ W(t,\theta)
\, W(t',\theta') \right] \, = \,
(t\wedge t')(\theta\wedge\theta'),
\qquad (t,\theta),(t',\theta')\in{\overline {\cal O}}.
\]
Let $\overline{W}_3:=(\overline{W}_3^i)_{i=1,2,3}$ be 
a ${\mathbb R}^3$-valued process, whose components
are three independent
copies of $W$, defined on $(\Omega,{\cal F}, {\mathbb P})$.
We denote by ${\cal F}_t$ the $\sigma$-field generated by
the random variables 
$(W(s,\theta):(s,\theta)\in[0,t]\times[0,1])$.

\vspace{.3cm}\noindent
We set for $\overline{x}\in H^3=L^2((0,1);{\mathbb R}^3)$:
\[
w_3(t,\theta):=\int_0^t\int_0^1 g_{t-s}(\theta,\theta') \, 
\overline{W}_3(ds,d\theta')
\]
\[
z_3(t,\theta):=e^{tA}\overline{x}(\theta)+w_3(t,\theta),
\qquad Z_3(t,\overline{x}):=z_3(t,\cdot).
\]
Then $z_3$ is the unique solution of the 
following ${\mathbb R}^3$-valued linear SPDE
with additive white-noise:
\begin{equation}\label{oud}
\left\{ \begin{array}{ll}
{\displaystyle
\frac{\partial z_3}{\partial t}=\frac 12
\frac{\partial^2 z_3}{\partial \theta^2}
 + \frac{\partial^2 \overline{W}_3}{\partial t\partial\theta} }
\\ \\
z_3(t,0)=z_3(t,1)=0
\\ \\
z_3(0,\theta)\, =\overline{x}(\theta) 
\end{array} \right.
\end{equation}
where $\overline{x}\in H^3$. The process
$z_3$ is also called the ${\mathbb R}^3$-valued random
string: see \cite{fun} and \cite{mutr}.
Recall that the law of $Z_3(t,\overline{x})$ is the
Gaussian measure ${\cal
N}(e^{tA}\overline{x},Q_t)$ on $H^3$, with
mean $e^{tA}\overline{x}$ and covariance operator
$Q_t:H^3\mapsto H^3$:
\begin{equation}\label{proqm_t}
Q_t\overline{h}(\theta) 
\, = \, \int_0^t \int_0^1 q_t(\theta,\theta')
\, \overline{h}(\theta') \, d\theta'\, ds, 
\end{equation}
for all $t\in[0,\infty]$, $\theta\in(0,1)$, $\overline{h}\in H^3$, where:
\begin{equation}\label{defq_t}
q_t(\theta,\theta') \, := \, \int_0^t g_{2s}(\theta,\theta') \, ds,
\quad t\in[0,\infty], \ \theta,\theta'\in(0,1).
\end{equation}
We denote by $(\overline{\beta}(\theta))_{\theta\in[0,1]}$
a 3-dimensional standard Brownian Bridge, and by 
$\mu_3$ the law of $\overline{\beta}$. Recall that
$\mu_3$ is equal to the Gaussian measure ${\cal N}(0,
Q_\infty)$ on $H^3$, $Q_\infty=(-2A)^{-1}$, and:
\begin{equation}\label{defq_inf}
q_\infty(\theta,\theta') = \theta\wedge\theta'-\theta\theta'.
\end{equation}
We set also for all $t\in[0,\infty)$, $\theta,\theta'\in(0,1)$:
\begin{equation}\label{defq^t}
q^t(\theta,\theta') \, := \, \int_t^\infty g_{2s}(\theta,\theta') \, ds
\, = \, q_\infty(\theta,\theta')-q_t(\theta,\theta').
\end{equation}
It is well
known that $Z_3$ is the Markov process associated with the
Dirichlet Form in $L^2(\mu_3)$:
\[
\Lambda^3(F,G):=\int_{H^3} \langle \overline{\nabla} F
, \overline{\nabla} G \rangle_{H^3} \,
d\mu_3, \qquad F,G\in W^{1,2}(\mu_3),
\]
where $\overline{\nabla} F:H^3\mapsto H^3$ is the gradient
of $F$ in $H^3$. For all $f:H^3\mapsto{\mathbb R}$ bounded and
Borel and for all $\overline{x}\in H^3$ we set:
\[
R_3(1)f(\overline{x}) \, := \, 
\int_0^\infty e^{-t} \,
{\mathbb E} \left[f(Z_3(t,\overline{x})) \right] \, dt.
\]

\vspace{.3cm}\noindent
The main result of this section is the following:
\begin{proposition}\label{preli}
$ $
\begin{itemize}
\item[1.] For all $\theta\in(0,1)$, $a\in{\mathbb R}^3$,
the function $U^{\theta,a}_3:H^3\mapsto
{\mathbb R}$:
\begin{equation}\label{defu}
U^{\theta,a}_3(\overline{x}) \, := \, \int_0^\infty e^{-t} \, 
\frac 1{(2\pi q_t(\theta,\theta))^{3/2}}
\exp\left(-\frac{|e^{tA}\overline{x}(\theta)-a|^2}{2q_t(\theta,\theta)}
\right) \, dt,
\end{equation}
is well defined and belongs to $C_b(H^3)\cap W^{1,2}(\mu_3)$.
If $(a_n,\theta_n)\to(a,\theta)\in{\mathbb R}^3\times (0,1)$, then:
\begin{equation}\label{unifa}
\lim_{n\to\infty}
\int_{H^3} \left[ \left|U^{\theta_n,a_n}_3-U^{\theta,a}_3\right|^2 +
\left\|\overline{\nabla} U^{\theta_n,a_n}_3-
\overline{\nabla} U^{\theta,a}_3\right\|^2 \right]\, d\mu_3
\, = \, 0.
\end{equation}
Moreover $(\theta^{3/2}(1-\theta)^{3/2}
U^{\theta,a}_3)_{\theta\in (0,1), \, a\in{\mathbb R}^3}$
is uniformly bounded, i.e.:
\begin{equation}\label{unibo}
\sup_{\theta\in(0,1), \, a\in{\mathbb R}^3}
\theta^{3/2}(1-\theta)^{3/2} \, \sup_{\overline{x}\in H^3}
\, U^{\theta,a}_3(\overline{x}) \, < \, \infty.
\end{equation}
\item[2.] Set $\overline{\gamma}^\theta(\overline{x}):=
|\overline{x}(\theta)|/\sqrt{\theta}$, $\overline{x}\in (C_0)^3$.
Then $\Gamma^\theta_3:=R_3(1)\overline{\gamma}^\theta$
converges to $\sqrt{8/\pi}$ in $W^{1,2}(\mu_3)$ as
$\theta\to 0$ or $\theta\to 1$.
\end{itemize}
\end{proposition}
{\bf Proof}. Let $\omega_3:=4\pi/3$. If
$\lambda\in{\mathbb R}$, we denote by $\lambda \cdot I$ the
linear application ${\mathbb R}^3\ni \alpha\mapsto \lambda\cdot
\alpha \in {\mathbb R}^3$. 

\vspace{.3cm}\noindent
{\bf Step 1}. Let $\overline{x}\in H^3$ be fixed. 
Notice that $z_3(t,\theta)$ has law
${\cal N}(e^{tA}\overline{x}(\theta),q_t(\theta,\theta)\cdot I)$,
where $q_t(\theta,\theta)$ is defined as in (\ref{defq_t}).
We denote by $(G_t(a,b):t,a,b>0)$ the fundamental solution
of the heat equation on $(0,+\infty)$ with homogeneous Dirichlet
boundary condition. By the reflection principle we have the
explicit representation:
\[
G_t(a,b) = \frac 1{\sqrt{2\pi t}}\exp\left(-\frac{(a-b)^2}{2t}\right)
\left(1-\exp\left(-\frac{2ab}t\right)\right).
\]
We set $\tau_{\theta}:=\inf\{t>0:\theta+B_t\in\{0,1\}\}$,
$\theta\in (0,1)$. Then we have:
\[
g_t(\theta,\theta')=G_t(\theta,\theta') - 
{\mathbb E}\left[1_{(t>\tau_{\theta'},\theta'+B_{\tau_{\theta'}}=1)}
G_{t-\tau_{\theta'}}(\theta,1)\right].
\]
Let $c_0:=1-\exp(-1)\in(0,1)$. Then for all $t>0$ and $a\geq 0$:
\[
\frac{c_0}{\sqrt{2\pi t}} \left(1\wedge\frac{2a^2}t\right)
\leq
G_t(a,a) \leq \frac 1{\sqrt{2\pi t}} \left(1\wedge\frac{2a^2}t\right).
\]
Let now $\theta\in[0,1/2]$. Then:
\begin{eqnarray*}
& &
\int_0^t
{\mathbb E}\left[1_{(2s>\tau_{\theta},\theta+B_{\tau_{\theta}}=1)}
G_{2s-\tau_{\theta}}(\theta,1)\right] \, ds
\\ \\ & = & {\mathbb E}\left[ \int_0^{(2t-\tau_{\theta})^+}
\frac{1_{(\theta+B_{\tau_{\theta}}=1)}}{2\sqrt{2\pi r}}
\exp\left(-\frac{(\theta-1)^2}{2r}\right)
\left(1-\exp\left(-\frac{2\theta}{r}\right)\right) dr
\right]
\\ \\
& \leq & {\mathbb E}\left[ 1_{(\theta+B_{\tau_{\theta}}=1)}\right]
\int_0^{2t} \frac 1{2\sqrt{2\pi r}}
\exp\left(-\frac{(\theta-1)^2}{2r}\right)
\left(1-\exp\left(-\frac{2\theta}{r}\right)\right)dr
\\ \\
& \leq & \frac{c_1}{\sqrt{\pi}} t\, \theta^2, \qquad
c_1 \, := \, \sup_{r>0} \frac 1{\sqrt{2r^3}}
\exp\left(-\frac{1}{8r}\right)\, < \, \infty.
\end{eqnarray*}
For all $t>0$ and $\theta\in[0,1/2]$ we obtain:
\begin{eqnarray*}
q_t(\theta,\theta) & \geq &
\int_0^t G_{2s}(\theta,\theta) \, ds - 
\frac {c_1}{\sqrt{\pi}} \, t\, \theta^2
\geq \int_0^t \frac{c_0}{2\sqrt{\pi s}} \left(1\wedge\frac{\theta^2}s\right)
\, ds - \frac {c_1}{\sqrt{\pi}} \, t\, \theta^2
\\ \\
& = & \frac{c_0}{\sqrt{\pi}}\left(
1_{(t\leq \theta^2)} \, {\sqrt t} \, + \, 1_{(\theta^2\leq t)}
\, \left(2\theta-\frac{\theta^2}{\sqrt t}\right)\right)
\, - t\, \frac {c_1}{\sqrt{\pi}} \, \theta^2
\\ \\
& \geq & \frac 1{\sqrt{\pi}}\left( c_0\, 
1_{(t\leq \theta^2)} \, {\sqrt t} \, + \, c_0\, 1_{(\theta^2\leq t)}
\, \theta \, -  \, c_1 \, t\, \theta^2 \right)
\end{eqnarray*}
Let $t_0:=(c_0/2c_1)\wedge(c_0/2c_1)^2$. If $t\geq t_0$, then
$q_t(\theta,\theta)\geq q_{t_0}(\theta,\theta)$. 
If $t\leq t_0$ then:
\begin{eqnarray*}
\frac \theta{q_t(\theta,\theta)} & \leq & \sqrt{\pi}\left(
\frac \theta{{\sqrt t}\,(c_0-c_1 {\sqrt t}\,\theta^2)} \, 1_{(t\leq \theta^2)}
 \, + \, \frac \theta{\theta(c_0-c_1 t\theta)} \, 1_{(\theta^2\leq t)}
\right)
\\ \\
&\leq & \frac {2\sqrt{\pi}}{c_0}
\left(\frac 1{\sqrt t} \, + \, 1 \right)
\end{eqnarray*}
By symmetry, we obtain that there exists $C_0>0$ such that
for all $\theta\in(0,1)$:
\begin{equation}\label{estq}
\left(\frac{\theta(1-\theta)}{q_t(\theta,\theta)} \right)^{3/2} \, \leq \,
C_0 \left(\frac 1{t^{3/4}} \wedge 1\right), \qquad t>0.
\end{equation}

\vspace{.3cm}\noindent
{\bf Step 2}. Fix $\theta\in(0,1)$. By (\ref{estq}), 
$U^{\theta,a}_3$ is well defined and in 
$C_b(H^3)$. Moreover for all $\overline{x} \in H^3$:
\begin{eqnarray*}
\theta^{3/2}(1-\theta)^{3/2} 
U^{\theta,a}_3(\overline{x}) & = &
\int_0^\infty e^{-t}
\left(\frac{\theta(1-\theta)}{2\pi q_t(\theta,\theta)}\right)^{3/2}
\exp\left(-\frac{|e^{tA}\overline{x}(\theta)-a|^2}{2q_t(\theta,\theta)}
\right) dt
\\ \\ & \leq & \frac{C_0}{(2\pi)^{3/2}}
\int_0^\infty e^{-t} \, \left(\frac 1{t^{3/4}} \wedge 1 \right)
\, dt \, < \, \infty,
\end{eqnarray*}
so that (\ref{unibo}) is proved. For all $\epsilon>0$ we set:
\[
\overline{f}^\epsilon(\overline{y}):=
\frac 1{\omega_3\epsilon^3} \,
1_{(|\overline{y}(\theta)-a|\leq \epsilon)},
\quad \overline{y}\in (C_0)^3.
\]
Let $\overline{x}\in H^3$. Then:
\begin{eqnarray}\label{ugu}
\nonumber & &
R_3(1)\overline{f}^\epsilon(\overline{x})
= \int_0^\infty e^{-t} \frac 1{\omega_3\epsilon^3}
{\mathbb P}\left(|z_3(t,\theta)+e^{tA}\overline{x}(\theta)-a|
\leq \epsilon \right) \, dt
\\ \nonumber \\ \nonumber 
& & = \, \int_0^\infty e^{-t} \, \frac 1{\omega_3\epsilon^3} \,
\int_{{\mathbb R}^3} 1_{(|\alpha|\leq \epsilon)} \, 
{\cal N}\left(e^{tA}\overline{x}(\theta)-a,q_t(\theta,\theta)\cdot I\right)
(d\alpha) \, dt
\\ \nonumber \\ \nonumber 
& & = \, \int_0^\infty dt \, e^{-t} \, \frac 1{\omega_3\epsilon^3} \,
\int_{(|\alpha|\leq \epsilon)} \frac 1{(2\pi q_t(\theta,\theta))^{3/2}}
\exp\left(-\frac{|\alpha-e^{tA}\overline{x}(\theta)+a|^2}
{2q_t(\theta,\theta)} \right) \, d\alpha 
\\ \nonumber \\ \nonumber 
& & = \frac 1{\omega_3\epsilon^3}
\int_{(|\alpha|\leq \epsilon)} \left[ \int_0^\infty e^{-t}
\frac 1{(2\pi q_t(\theta,\theta))^{3/2}}
\exp\left(-\frac{|\alpha-e^{tA}\overline{x}(\theta)+a|^2}{2q_t(\theta,\theta)}
\right) dt \right] d\alpha
\\ \nonumber \\
& & = \, \frac 1{\omega_3\epsilon^3}
\int_{(|\alpha|\leq \epsilon)}
U^{\theta,a+\alpha}_3(\overline{x}) \, d\alpha.
\end{eqnarray}
By (\ref{estq}) and the Dominated Convergence Theorem, we have that
for all $(\theta,a)\in(0,1)\times{\mathbb R}^3$:
\begin{equation}\label{ugu1}
\lim_{\epsilon\to 0} R_3(1)\overline{f}^\epsilon(\overline{x}) \,
= \, U^{\theta,a}_3(\overline{x}), \qquad
\forall \, \overline{x}\in H^3,
\end{equation}
uniformly for $\overline{x}$ in bounded sets of $H^3$, and
by (\ref{estq}):
\begin{equation}\label{ugu2}
|R_3(1)\overline{f}^\epsilon(\overline{x})| \, \leq \,
\int_0^\infty e^{-t}
\frac 1{(2\pi q_t(\theta,\theta))^{3/2}} \, dt \, < \, \infty.
\end{equation}

\vspace{.3cm}\noindent
{\bf Step 3}. 
Notice that by the Dominated Convergence Theorem the map
${\mathbb R}^3\times (0,1)\ni(a,\theta) \mapsto
U^{\theta,a}_3 \in L^2(\mu_3)$ is continuous.
We want to prove now that $U^{\theta,a}_3$ is in
$W^{1,2}(\mu_3)$: to this aim we shall prove that
$R_3(1)\overline{f}^\epsilon$
converges to $U^{\theta,a}_3$ in $W^{1,2}(\mu_3)$.
We define for all $\overline{x}\in (C_0)^3$,
$a\in{\mathbb R}^3\backslash\{\overline{x}(\theta)\}$: 
\[
{\cal U}^{\theta,a}_{\, \overline{h}}(\overline{x})
\, := \, - \, \int_0^\infty e^{-t} \, 
\frac{e^{tA}\overline{h}(\theta)}
{(2\pi)^{3/2}(q_t(\theta,\theta))^{2}} \, 
\psi\left((e^{tA}\overline{x}(\theta)-a)/
{\sqrt{q_t(\theta,\theta)}}\right) dt,
\]
\begin{equation}\label{defpsi}
{\rm where}: \qquad 
\psi:{\mathbb R}^3\mapsto{\mathbb R}^3, \quad
\psi(a) \, := \, a \, \exp\left(-\frac{|a|^2}2\right).
\end{equation}
Recall now that for all $\overline{h}\in H^3$ and $\theta\in(0,1)$:
\begin{eqnarray}\label{estgr}
& & \nonumber
|e^{tA}\overline{h}(\theta)| \, \leq \, e^{tA}|\overline{h}|(\theta)
\, = \int_0^1 |\overline{h}(\theta')| \, g_t(\theta,\theta') \, d\theta'
\, \leq \int_0^1
|\overline{h}(\theta')| \, G_t(\theta,\theta') \, d\theta'
\\ \nonumber \\ & &
\leq \, \frac{1\wedge(2\theta/t)}{\sqrt{2\pi t}}\int_0^1
|\overline{h}(\theta')| \, 
\exp\left( -\frac{|\theta-\theta'|^2}{2t}\right) d\theta'
\, \leq \,
\frac{1\wedge(2\theta/t)}{t^{1/4}} \,
\|\overline{h}\|,
\end{eqnarray}
so that:
\[
\left|\sup_{\|\overline{h}\|=1} 
{\cal U}^{\theta,a}_{\, \overline{h}} (\overline{x})\right|
\, \leq \, \int_0^\infty \frac{e^{-t}}
{(q_t(\theta,\theta))^{2}t^{1/4}}
\left| \psi\left((e^{tA}\overline{x}(\theta)-a)/
{\sqrt{q_t(\theta,\theta)}}\right) \right| \, dt.
\]
Since $\overline{\beta}$ has law $\mu_3={\cal N}(0,Q_\infty)$,
then $e^{tA}\overline{\beta}$ has law ${\cal
N}(0,e^{tA}Q_\infty e^{tA})={\cal N}(0,Q_\infty-Q_t)$.
Then, by (\ref{estq}), and since $|\psi|\leq 1$:
\begin{eqnarray*}
& &
\left({\mathbb E}\left[
\left|\sup_{\|\overline{h}\|=1}
{\cal U}^{\theta,a}_{\, \overline{h}}
\left(\overline{\beta}\right)\right|^2 \right] \right)^{1/2}
\, \leq \, \int_0^\infty \frac{e^{-t}}
{(q_t(\theta,\theta))^{2}t^{1/4}} \, \cdot
\\ \\ & &
\quad \cdot \,
\left( \int_{{\mathbb R}^3}\left|
\psi\left(\alpha'/
{\sqrt{q_t(\theta,\theta)}}\right) \right|^2 \,
{\cal N}(a,q^t(\theta,\theta)\cdot I)(d\alpha')\right)^{1/2}\, dt 
\\ \\ & &
\leq \, \int_0^\infty \frac{e^{-t}}
{(q_t(\theta,\theta))^{2}t^{1/4}} \,
\left\{1 \wedge \left[\left(\frac{q_t(\theta,\theta)}
{2\pi q^t(\theta,\theta)}\right)^{3/4} \,
\|\psi\|_{L^2({\mathbb R}^3)}\right] \right\}
\, dt 
\\ \\ & &
\leq \, \frac 1{(q^1(\theta,\theta))^{3/4}} \, 
\int_0^1 \frac 1{(q_t(\theta,\theta))^{5/4}t^{1/4}}
 \, dt \, + \, \frac 1{(q_1(\theta,\theta))^2} \,
\int_1^\infty e^{-t} \, dt \, < \, \infty.
\end{eqnarray*}
Therefore, setting for $\mu_3$-a.e. $\overline{x}$:
\[
{\cal U}^{\theta,a}
\, := \, - \, \int_0^\infty e^{-t} \, 
\frac{g_t(\theta,\cdot)}
{(2\pi)^{3/2}(q_t(\theta,\theta))^{2}} \, 
\psi\left((e^{tA}\overline{x}(\theta)-a)/
{\sqrt{q_t(\theta,\theta)}}\right) dt,
\]
we have that ${\cal U}^{\theta,a}
\in L^2(H^3,\mu_3\, ;H^3)$, and
$\langle{\cal U}^{\theta,a},\overline{h}
\rangle \, = \, {\cal U}^{\theta,a}_{\, \overline{h}}$
in $L^2(\mu_3)$, for all  $\overline{h}\in H^3$.
Arguing analogously we have:
\begin{eqnarray}\label{limnab}
& & \nonumber
\left({\mathbb E}\left[\left\| 
{\cal U}^{\theta,a+\alpha} \left(\overline{\beta}\right)-
{\cal U}^{\theta,a} \left(\overline{\beta}\right)
\right\|^2 \right] \right)^{1/2} 
\, \leq \, \int_0^\infty \frac{e^{-t}}
{(q_t(\theta,\theta))^{2}t^{1/4}} \cdot
\\ \nonumber \\ & & \cdot 
\left\{1 \wedge \left[\left(\frac{q_t(\theta,\theta)}
{2\pi q^t(\theta,\theta)}\right)^{3/4}
\left\|\psi\left(\cdot+\alpha/{\sqrt{q_t(\theta,\theta)}}
\right)-\psi\right\|_{L^2({\mathbb R}^3)}\right] \right\}
dt \qquad
\end{eqnarray}
which tends to $0$ as $\alpha\to 0$. 
Therefore we can differentiate under the integral sign in
(\ref{ugu}) and obtain:
\[
\overline{\nabla} R_3(1)\overline{f}^\epsilon
\, = \, 
\frac 1{\omega_3\epsilon^3} 
\int_{(|\alpha|\leq \epsilon)} 
{\cal U}^{\theta,a+\alpha} \, d\alpha, \qquad {\rm in} \
L^2(H^3,\mu_3\, ;H^3).
\]
Therefore by (\ref{limnab}):
\begin{eqnarray*}
& &
\int_{H^3}
\left\| \overline{\nabla} R_3(1)\overline{f}^\epsilon-
{\cal U}^{\theta,a} \right\|^2  d\mu_3
\\ \\ & & \leq \,
\frac 1{\omega_3\epsilon^3} \int_{(|\alpha|\leq \epsilon)} 
\int_{H^3} \left\| {\cal U}^{\theta,a+\alpha}
- {\cal U}^{\theta,a} \right\|^2
d\mu_3 \, d\alpha  \, \to \, 0
\end{eqnarray*}
as $\epsilon\to 0$. Therefore, 
$R_3(1)\overline{f}^\epsilon$ converges to
$U^{\theta,a}_3$ in $L^2(\mu_3)$ and 
$\overline{\nabla} R_3(1)\overline{f}^\epsilon$ converges to 
${\cal U}^{\theta,a}$ in $L^2(H^3,\mu_3\, ;H^3)$ 
as $\epsilon\to 0$. Since $W^{1,2}(\mu_3)$
is complete, then $U^{\theta,a}_3\in W^{1,2}(\mu_3)$, 
$\overline{\nabla} U^{\theta,a}_3 \, = \, {\cal U}^{\theta,a}$ in
$L^2(H^3,\mu_3\, ;H^3)$ and
$R_3(1)\overline{f}^\epsilon$ converges to 
$U^{\theta,a}_3$ in $W^{1,2}(\mu_3)$ as $\epsilon\to 0$.  
Moreover, by (\ref{limnab}), (\ref{unifa}) is proved.

\vspace{.3cm}\noindent
{\bf Step 4}. We prove now the last assertion. By symmetry,
it is enough to consider the case $\theta\to 0$.
Recall that $\overline{\gamma}^\theta(\overline{x})=
|\overline{x}(\theta)|/\sqrt{\theta}$, $\overline{x}\in (C_0)^3$. Then:
\begin{eqnarray*}
\Gamma^\theta_3(\overline{x}) & := & R_3(1)
\overline{\gamma}^\theta(\overline{x})
\, = \, \frac 1{\sqrt{\theta}}
\int_0^\infty e^{-t} \int_{{\mathbb R}^3}
|\alpha|\, {\cal N}\left( e^{tA}\overline{x}
(\theta), q_t(\theta,\theta) \cdot I \right)(d\alpha) \, dt
\\ \\ & = & 
\int_0^\infty e^{-t} \sqrt{\frac{q_t(\theta,\theta)}\theta}
\int_{{\mathbb R}^3} |\alpha|\, {\cal N}\left( 
e^{tA}\overline{x}(\theta) / \sqrt{q_t(\theta,\theta)},
I \right)(d\alpha) \, dt,
\end{eqnarray*}
and:
\begin{eqnarray*}
& &
\overline{\nabla}\Gamma^\theta_3(\overline{x}) \, = \, - \,
\int_0^\infty e^{-t} \, \sqrt{\frac{q_t(\theta,\theta)}\theta} \, 
g_t(\theta,\cdot) \, \cdot
\\ \\ & & \quad \cdot \,
\int_{{\mathbb R}^3} |\alpha|\, \left(\alpha-
e^{tA}\overline{x}(\theta) / \sqrt{q_t(\theta,\theta)} \right)
{\cal N}\left( 
e^{tA}\overline{x}(\theta) / \sqrt{q_t(\theta,\theta)},
I \right)(d\alpha) \, dt,
\end{eqnarray*}
for all $\overline{x}\in (C_0)^3$. By (\ref{estgr}) and
by Schwartz's inequality:
\[
\|\overline{\nabla}\Gamma^\theta_3(\overline{x})\| \, \leq \,
\sqrt{3} \int_0^\infty e^{-t} \, \frac 1{t^{1/4}} \,
\left(1 \wedge\frac{2\theta}t\right) \, 
\sqrt{3 + \frac{|e^{tA}\overline{x}(\theta)|^2}
{q_t(\theta,\theta)} } \ dt.
\]
By the sub-additivity of the square-root, by (\ref{estq})
and since $q^t(\theta,\theta)\leq \theta(1-\theta)$:
%
\begin{eqnarray*}
& &
\left[{\mathbb E}\left(\left\|\overline{\nabla}\Gamma^\theta_3
\left(\overline{\beta}\right)\right\|^2\right)\right]^{1/2} \, \leq \,
\sqrt{3} \int_0^\infty  \frac{e^{-t}}{t^{1/4}} \,
\left(1 \wedge\frac{2\theta}t\right) \, \left(
\sqrt{3} + \sqrt{\frac{q^t(\theta,\theta)}
{q_t(\theta,\theta)}} \right) dt
\\ \\ & & \leq \, \sqrt{3} \int_0^\infty e^{-t} \,
\left(1 \wedge\frac{2\theta}t\right) 
\left( \frac{\sqrt 3}{t^{1/4}} \, + \,
\frac{(C_0)^{1/3}}{t^{3/4}}\right) dt \, \to \, 0
\end{eqnarray*}
as $\theta\to 0$. Since $\mu_3$ is invariant for $z_3$, we
have
\[
{\mathbb E}\left(
\Gamma^\theta_3\left(\overline{\beta}\right) \right) \, = \,
\frac 1{\sqrt{\theta}} \, 
{\mathbb E}\left(|\overline{\beta}(\theta)|\right) \, =: \,
c_\theta.
\]
By the Poincar\'e inequality for $\Lambda^3$, see \cite{dpz3},
there exists $C>0$ such that:
\[
{\mathbb E}\left(\Gamma^\theta_3
\left(\overline{\beta}\right)-c_\theta\right)^2 \, \leq \,
\frac 1 C \, {\mathbb E}\left(\left\|\overline{\nabla}\Gamma^\theta_3
\left(\overline{\beta}\right)\right\|^2\right) \, \to \, 0,
\]
as $\theta\to 0$, and since:
\begin{eqnarray*}
& &
c_\theta \, = \, \frac{4\pi}{\sqrt\theta} \int_0^\infty 
\frac 1{(2\pi\theta(1-\theta))^{3/2}} 
\, r^3 \, \exp\left(-\frac{r^2}{2\theta(1-\theta)}\right) \, dr
\\ \\ & & = \, \sqrt{1-\theta} \, 
\sqrt{\frac 2{\pi}}
\int_0^\infty r^3 \exp\left(-\frac{r^2}{2}\right)\, dr \, = \,
\sqrt{1-\theta} \, \sqrt{\frac 8{\pi}} 
\, \to \, \sqrt{\frac 8\pi}
\end{eqnarray*}
we obtain that $\Gamma^\theta_3$ converges to
$\sqrt{8/\pi}$ in $W^{1,2}(\mu_3)$ as $\theta\to 0$.
\quad $\square$

\section{Occupation densities}

Following \cite{nupa}, we set the:
\begin{definition}\label{d1}
A pair $(u,\eta)$ is said to be a
solution of equation (\ref{1}) 
with initial value $x\in K_0 \cap C_0$, if:
\begin{itemize}
\item[(i)] $\{u(t,\theta):(t,\theta)\in {\overline{\cal O}}\}$ is a continuous
and adapted process, i.e. $u(t,\theta)$ is ${\cal F}_t$-measurable
for all $(t,\theta)\in {\overline{\cal O}}$, a.s. 
$u(\cdot,\cdot)$ is continuous on
${\overline{\cal O}}$, $u(t,\cdot)\in K_0\cap C_0$ 
for all $t\geq 0$, and $u(0,\cdot)=x$.
\item[(ii)] $\eta(dt,d\theta)$ is a random positive measure on 
${\cal O}$ such that $\eta([0,T]\times[\delta,1-\delta])<+\infty$
for all $T,\delta>0$, and $\eta$ is adapted, i.e. $\eta(B)$ is
${\cal F}_t$-measurable for every Borel set $B\subset 
[0,t]\times(0,1)$.
\item[(iii)] For all $t\geq 0$ and $h \in D(A)$
\begin{eqnarray}
\label{weeq}
\nonumber & &
\langle u(t,\,\cdot\,),h \rangle \, -
\, \langle x,h \rangle \, - \, \int_0^t 
\langle u(s,\,\cdot\,),Ah \rangle \, ds  
\\ \\ \nonumber & & = \, - \,
\int_0^1 h'(\theta) \, W(t,\theta) \, d\theta
+ \, \int_0^t\int_0^1 h(\theta)\, 
\eta(ds,d\theta), \qquad {\rm a.s.}.
\end{eqnarray}
\item[(iv)] $\int_{\cal O} u \ d\eta \ =\ 0$.
\end{itemize}
\end{definition}
In \cite{nupa}, existence and uniqueness solutions of
equation (\ref{1}) were proved.

\vspace{.3cm}\noindent
We denote by $(e(\theta))_{\theta\in[0,1]}$
the 3-Bessel Bridge between $0$ and $0$,
see \cite{reyo}, and by $\nu$ the law on $K_0$ of $e$.
We recall the following result, proved in \cite{za00}.
\begin{theorem}\label{zumpa}
Let $\Phi_3:H^3=L^2(0,1;{\mathbb R}^3)\mapsto K_0$, 
\ $\Phi_3(y)(\theta):=|y(\theta)|_{{\mathbb R}^3}$. 
\begin{itemize}
\item[1.] The process $u$ is a Strong-Feller Markov
process properly associated with the symmetric Dirichlet Form
${\cal E}$ in $L^2(\nu)$:
\[
\frac 12 \, \int_{K_0} \langle \nabla \varphi,
\nabla \psi \rangle \, d\nu, \qquad
\varphi,\psi \, \in \, W^{1,2}(\nu).
\]
\item[2.] The Dirichlet Form ${\cal E}$
is the image of $\Lambda^3$ under the map $\Phi_3$, i.e.
$\nu$ is the image $\mu_3$ under $\Phi_3$ and:
\[
W^{1,2}(\nu)=\{\varphi\in L^2(\nu)
: \varphi\circ\Phi_3\in W^{1,2}(\mu_3)\},
\]
\begin{equation}\label{po3z2}
{\cal E}(\varphi,\psi) \, = \, \Lambda^3
(\varphi\circ\Phi_3,\psi\circ\Phi_3) \quad
\forall \, \varphi,\psi \in W^{1,2}(\nu).
\end{equation}
\end{itemize}
\end{theorem}
We refer to
\cite{fot} and \cite{maro} for all basic definitions in the
Theory of Dirichlet Forms. Notice that by point 1 in Theorem
\ref{zumpa} and by Theorem IV.5.1 in \cite{maro}, the Dirichlet
Form ${\cal E}$ is quasi-regular. In particular, by the transfer
method stated in VI.2 of \cite{maro} we can apply several results
of \cite{fot} in our setting.

\vspace{.3cm}\noindent
We recall the definition of an
Additive Functional of the Markov process $u$. We denote by
$({\mathbb P}_x:x\in K_0)$ the family the of laws of $u$
on $E:=C([0,\infty);K_0)$ and the coordinate process on $K_0$
by: $X_t:E\mapsto K_0$, $t\geq 0$,
$X_t(e):=e(t)$.
By a Positive Continuous Additive Functional (PCAF) in the strict
sense of $u$, we mean a
family of functions $A_t:E\mapsto {\mathbb R}^+$,
$t\geq 0$, such that:
\begin{itemize}
\item[(A.1)] $(A_t)_{t\geq 0}$ is adapted to the minimum
admissible filtration $({\cal N}_t)_{t\geq 0}$ of $u$, see Appendix
A.2 in \cite{fot}.
\item[(A.2)] There exists a set $\Lambda\in{\cal N}_\infty$ such that
${\mathbb P}_x(\Lambda)=1$ for all $x\in K_0$,
$\theta_t(\Lambda)\subseteq \Lambda$ for all $t\geq 0$, and for all
$\omega\in \Lambda$: $t\mapsto A_t(\omega)$ is 
continuous non-decreasing, $A_0(\omega)= 0$ and for all $t,s\geq 0$:
\begin{equation}\label{add1}
A_{t+s}(\omega) \, = \, A_s(\omega)+A_t(\theta_s\omega),
\end{equation}
where $(\theta_s)_{s\geq 0}$ is the time-translation semigroup on
$E$. 
\end{itemize}
Two PCAFs in the strict sense 
$A^1$ and $A^2$ are said to be equivalent if
\[
{\mathbb P}_x \left( A^1_t=A^2_t \right) \, = \, 1, 
\quad \forall t>0, \ \forall x\in K_0.
\] 
If $A$ is a linear combination of PCAFs in the strict sense of 
$u$, then the Revuz-measure of $A$ 
is a Borel signed measure $m$ on $K_0$ such that:
\[
\int_{K_0} \varphi\, dm \, = \, 
\int_{K_0} {\mathbb E}_x\left[
\int_0^1 \varphi(X_t)\, dA_t \right]\,
\nu(dx), \quad \forall \varphi\in C_b(K_0).
\]
Moreover $U\in D({\cal E})$ is the 1-potential of a
PCAF $A$ in the strict sense with Revuz-measure $m$, if:
\[
{\cal E}_1(U,\varphi) \, = \, \int_{K_0} \varphi \, dm,
\quad  \forall \ \varphi\in D({\cal E}) \cap C_b(K_0),
\]
where ${\cal E}_1:={\cal E}+(\cdot,\cdot)_{L^2(\nu)}$.
We introduce the following notion of convergence of Positive
Continuous Additive Functionals in the strict sense of $X$.
\begin{definition}\label{def1}
Let $(A_n(t))_{t\geq 0}$, $n\in{\mathbb N}\cup\{\infty\}$,
be a sequence of PCAF's in the strict sense of $u$. We say
that $A_n$ converges to $A_\infty$, if:
\begin{itemize}
\item[1.] For all $\epsilon>0$ and for all $x\in K_0\cap C_0$:
\begin{equation}\label{defe1}
\lim_{n\to\infty} A_n(t+\epsilon)-A_n(\epsilon) \, = \,
A_\infty(t+\epsilon)-A_\infty(\epsilon),
\end{equation}
uniformly for $t$ in compact sets of $[0,\infty)$,
${\mathbb P}_x$-almost surely. 
\item[2.] For ${\cal E}$-q.e. $x\in K_0\cap C_0$:
\begin{equation}\label{defe2}
\lim_{n\to\infty} A_n(t) \, = \, A_\infty(t),
\end{equation}
uniformly for $t$ in compact sets of $[0,\infty)$,
${\mathbb P}_x$-almost surely. 
\end{itemize}
\end{definition}
\begin{lemma}\label{lemdef}
Let $(A_n(t))_{t\geq 0}$, $n\in{\mathbb N}\cup\{\infty\}$,
be a sequence of PCAF's in the strict sense of $X$, and let
$U_n$ be the 1-potential of $A_n$, $n\in{\mathbb N}\cup\{\infty\}$.
If $U_n\to U_\infty$ in $D({\cal E})$, then $A_n$ converges
to $A_\infty$ in the sense of Definition \ref{def1}.
\end{lemma}
{\bf Proof.} Since $U_n\to U_\infty$ in $D({\cal E})$, by 
Corollary 5.2.1 in \cite{fot}, we have point 2 of 
Definition \ref{def1}, i.e. there exists an ${\cal E}$-exceptional set
$V$ such that (\ref{defe2}) holds for all $x\in K_0\backslash V$. By the 
Strong Feller property of $X$, ${\mathbb P}_x$-a.s.
$X_t\in E\backslash V$,
for all $t>0$ and for all $x\in K_0$, and by the 
additivity property (\ref{defe1}) holds for all $x\in K_0$.
\quad $\square$
\begin{remark}\label{meas}
{\rm
We recall that if $(A,{\cal A})$ is a 
measurable space,
$(\Omega,{\cal F},{\mathbb P})$ a probability space and
$X_n$ is a sequence of ${\cal A}\otimes{\cal F}$-measurable
random variables, such that $X_n(a,\cdot)$ converges in
probability for every $a\in A$, then there exists a 
${\cal A}\otimes{\cal F}$-measurable random variable
$X$, such that $X(a,\cdot)$ is the limit in probability
of $X_n(a,\cdot)$ for every $a\in A$.
}
\end{remark}
We can now state the main result of this section:
\begin{theorem}\label{occuref}
Let $\theta\in(0,1)$, $a\geq 0$.
\begin{itemize}
\item[1.] For all $(\theta,a)\in(0,1)\times
[0,\infty)$, there exists a  PCAF in the strict sense
of $u$, $(l^a(t,\theta))_{t\geq 0}$, such that
$(l^a(\cdot,\theta))_{\theta\in(0,1), a\in[0,\infty)}$
is continuous in the sense of Definition \ref{def1}
and jointly measurable, and such that for all $a\geq 0$:
\[
l^a(t,\theta)
\, = \, \lim_{\epsilon\downarrow 0} \,
\frac 1{\epsilon} 
\int_0^t 1_{[a,a+\epsilon]}(u(s,\theta)) \, ds, \qquad t\geq 0,
\]
in the sense of Definition \ref{def1}.
\item[2.] The Revuz measure of $l^a(\cdot,\theta)$ is:
\[
\sqrt{\frac 2{\pi\theta^3(1-\theta)^3}} \ a^2
\, \exp\left(-\frac{a^2}{2\theta(1-\theta)}\right) 
\, \nu(\, dx \,| \, x(\theta)=a), \qquad a\geq 0,
\]
and in particular $l^0(\cdot,\theta)\equiv 0$. Moreover,
$l^a(\cdot,\theta)$ increases only on $\{t:u(t,\theta)=a\}$.
\item[3.] The following occupation times formula holds for all
$\theta\in(0,1)$:
\begin{equation}\label{otfr}
\int_0^t F(u(s,\theta)) \, ds \, = \, 
\int_0^\infty F(a) \, l^a(t,\theta) 
\, da, \quad F\in B_b({\mathbb R}), \ t \geq 0.
\end{equation}
\end{itemize}
\end{theorem}
For an overview on existence of occupation densities see \cite{geho}.

\vspace{.3cm}\noindent
We set $\Lambda_1^3:=\Lambda^3+(\cdot,\cdot)_{L^2(\mu_3)}$,
${\cal E}_1:={\cal E}+(\cdot,\cdot)_{L^2(\nu)}$.
For all $f:H\mapsto{\mathbb R}$ bounded and
Borel and for all $x\in K_0\cap C_0$ we introduce the
1-resolvent of $u$:
\[
R(1)f(x) = 
\int_0^\infty e^{-t} \, {\mathbb E}_x \left[f(X_t) \right]
\, dt,
\]
where ${\mathbb E}_x$ denotes the expectation w.r.t. the
law of the solution $u$ of (\ref{1}) with initial value $x$.
The next Lemma gives the projection principle from the
Dirichlet space $W^{1,2}(\mu_3)$, 
associated with the Gaussian process $z_3$,
to the Dirichlet space $W^{1,2}(\nu)$
of the solution $u$ of the SPDE
with reflection (\ref{1}).
\begin{lemma}\label{projector}
There exists a unique bounded linear operator $\Pi:
W^{1,2}(\mu_3)\mapsto W^{1,2}(\nu)$, such that for all
$F,G\in W^{1,2}(\mu_3)$ and $f\in W^{1,2}(\nu)$:
\begin{equation}\label{pro1}
\Lambda^3_1(F, f\circ\Phi_3) \, = \, 
{\cal E}_1(\Pi F,f),
\end{equation}
\begin{equation}\label{pro2}
\Lambda^3_1((\Pi F)\circ\Phi_3, G) \, = \, 
\Lambda^3_1(F,(\Pi G)\circ\Phi_3).
\end{equation}
In particular, we have that for all $\varphi\in L^2(\nu)$
and $F\in W^{1,2}(\mu_3)$:
\begin{equation}\label{pro3}
R(1)\varphi
\, = \, \Pi\left(R_3(1) \left[\varphi\circ\Phi_3\right] \right),
\end{equation}
\begin{equation}\label{estipro}
\|\Pi F\|_{{\cal E}_1}
\, \leq \,
\|F\|_{\Lambda^3_1}.
\end{equation}
Finally, $\Pi$ is Markovian, i.e. $\Pi 1=1$ and:
\begin{equation}\label{markov}
F\in W^{1,2}(\mu_3), \quad 0\leq F \leq 1 \ \Longrightarrow \
0\leq \Pi F\leq 1.
\end{equation}
\end{lemma}
{\bf Proof.} Let ${\cal D}:=\{\varphi\circ\Phi_3:\varphi\in
W^{1,2}(\nu)\}\subset W^{1,2}(\mu_3)$. Let $W^{1,2}(\mu_3)$
be endowed with the scalar product $\Lambda^3_1$: then, by
(\ref{po3z2}), ${\cal D}$ is a closed subspace of
$W^{1,2}(\mu_3)$. Therefore there exists 
a unique bounded linear projector $\hat{\Pi}:
W^{1,2}(\mu_3)\mapsto{\cal D}$, symmetric with respect to
the scalar product $\Lambda^3_1$. For all $F\in W^{1,2}(\mu_3)$
we set $\Pi F:=f$ where $f$ is the 
unique element of $W^{1,2}(\nu)$ such that
$f\circ\Phi_3 = \hat{\Pi}F$. Then (\ref{pro1}) and
(\ref{pro2}) are satisfied by construction. Let now
$\varphi,\psi\in W^{1,2}(\nu)$. Then by (\ref{po3z2}):
\begin{eqnarray*}
& &
{\cal E}_1(R(1)\varphi,\psi) \, = \,
\int_{K_0} 
\varphi\,\psi\, d\nu \, = \, \int_{H^3} (\varphi\circ\Phi_3)
\, (\psi\circ\Phi_3) \, d\mu_3
\\ \\ & & = \,
\Lambda^3_1(R_3(1)\left[\varphi\circ\Phi_3\right],
\psi\circ\Phi_3) \, = \,
\Lambda^3_1(\hat{\Pi}R_3(1)\left[\varphi\circ\Phi_3\right],
\psi\circ\Phi_3)
\\ \\ & & = \,
{\cal E}_1(\Pi R_3(1)\left[\varphi\circ\Phi_3\right],
\psi),
\end{eqnarray*}
which implies (\ref{pro3}). Then, since $\hat{\Pi}$ is a
symmetric projector:
\[
\|\Pi F\|_{{\cal E}_1} \, = \, \|\hat{\Pi}F\|_{\Lambda_1^3}
\, \leq \, \|F\|_{\Lambda_1^3},
\]
so that (\ref{estipro}) is proved. Notice now that $1\in {\cal
D}$, so that obviously $\Pi 1=1$. Moreover, recall that
$\hat{\Pi}F$ is characterized by the property: 
\[
\hat{\Pi}F\in{\cal D}, \quad
\Lambda_1^3(F-\hat{\Pi}F,G) \, = \, 0, \quad \forall 
G\in{\cal D}.
\]
Let $F\in W^{1,2}(\mu_3)$ such that $F\geq 0$. Since
${\cal E}$ is a Dirichlet Form, then $(\hat{\Pi}F)^-
:=(-\hat{\Pi}F)\vee 0$ still
belongs to ${\cal D}$, and since $\Lambda_1^3$
is a Dirichlet Form:
\[
0 \, = \, \Lambda_1^3(F-\hat{\Pi}F,(\hat{\Pi}F)^-) \,
= \, \Lambda_1^3(F,(\hat{\Pi}F)^-) \,
+ \, \|(\hat{\Pi}F)^-\|_{\Lambda_1^3}^2 \, \geq \,
\|(\hat{\Pi}F)^-\|_{\Lambda_1^3}^2,
\]
so that $\hat{\Pi}F\geq 0$, and (\ref{markov}) 
follows. \quad $\square$

\vspace{.3cm}
\noindent
{\bf Proof of Theorem \ref{occuref}.} 
Let $a\geq 0$. For all $\epsilon>0$ we set:
\[
f^\epsilon(y):=
\frac 1{\epsilon} \,
1_{[a,a+\epsilon]}(y(\theta)),
\quad y\in K_0\cap C_0.
\]
By Lemma \ref{projector}, we have that:
\begin{eqnarray*}
R(1)f^\epsilon & = &
\Pi \left( R_3(1)\left[f^\epsilon\circ\Phi_3\right]
\right)
\, = \, \frac 1\epsilon
\int_{(a\leq|\alpha|\leq a+\epsilon)}
\Pi \left( U^{\theta,a+\alpha}_3 \right)
\, d\alpha
\\ \\ & = & 
\frac 1\epsilon \int_a^{a+\epsilon} r^2 \, dr
\int_{{\mathbb S}^{2}} \Pi \left( U_3^{\theta,r\cdot
n}\right)  \, {\cal H}^{2}(dn),
\end{eqnarray*}
where ${\cal H}^{2}$ is the $2$-dimensional Hausdorff
measure, $U_3^{\theta,a\cdot n}$ is the 1-potential
in $W^{1,2}(\mu_3)$ defined by (\ref{defu}) 
and $\Pi$ is the operator defined in Lemma
\ref{projector}. By (\ref{unifa}) above,
the map $r\mapsto U_3^{\theta,r\cdot n}\in W^{1,2}(\mu_3)$
is continuous. Let $U^{\theta,a}\in W^{1,2}(\nu)$ be 
defined by:
\[
U^{\theta,a} \, := \, a^2 \, 
\int_{{\mathbb S}^{2}}
\Pi\left( U_3^{\theta,a\cdot n}\right) \,
{\cal H}^{2}(dn), \qquad a\geq 0.
\]
By (\ref{estipro}) we have that $R(1)f^\epsilon$ converges 
to $U^{\theta,a}$ in $W^{1,2}(\nu)$ as $\epsilon\to 0$.
For all $\epsilon>0$ and $\varphi\in W^{1,2}(\nu)\cap C_b(K_0)$
we have:
\[
{\cal E}_1(R(1)f^\epsilon,\varphi) \, = \, 
\int_{K_0} f^\epsilon \, \varphi \, d\nu \, = \, 
\frac 1\epsilon \, {\mathbb E}\left[ \varphi(e) \, 1_{[a,
a+\epsilon]}(e(\theta))\right],
\]
where the law of $e$ is $\nu$ and
${\cal E}_1={\cal E}+(\cdot,\cdot)_{L^2(\nu)}$.
Letting $\epsilon\to 0$ we get:
\begin{eqnarray}\label{smooth}
\nonumber
& &
{\cal E}_1(U^{\theta,a},\varphi) \, = \, \lim_{\epsilon\to 0}
\frac 1\epsilon \, {\mathbb E}\left[ \varphi(e) \, 1_{[a,
a+\epsilon]}(e(\theta))\right]
\\ \nonumber \\ & & = \,
\sqrt{\frac 2{\pi\theta^3(1-\theta)^3}} \ a^2
\, \exp\left(-\frac{a^2}{2\theta(1-\theta)}\right) 
\, {\mathbb E}\left[ \varphi(e) \, | \, e(\theta)=a\right].
\end{eqnarray}
By Lemma \ref{projector}, $\Pi$ is a Markovian operator
and by (\ref{unibo}) in Proposition \ref{preli} the family
$(U_3^{\theta,a\cdot n}:n\in{\mathbb S}^{2})$ is uniformly
bounded in the supremum-norm. Therefore, $U^{\theta,a}$ is
bounded, and by (\ref{smooth}) $U^{\theta,a}$ is the
1-potential of a non-negative finite measure. By Theorem 5.1.6
in \cite{fot}, there exists a PCAF $(l^a(t,\theta))_{t\geq
0}$ in the strict sense of $u$, with 1-potential equal to
$U^{\theta,a}$ and with Revuz-measure given by
(\ref{smooth}). 
Notice now that $R(1)f^\epsilon$ is the 1-potential of the 
following PCAF in the strict sense of $u$:
\[
t \, \mapsto \, \frac 1{\epsilon} 
\int_0^t 1_{[a,a+\epsilon]}(u(s,\theta)) \, ds, \qquad t\geq 0.
\]
Therefore, points 1 and 2 of Theorem \ref{occuref} are proved 
by (\ref{unifa}), Lemma \ref{lemdef} and Remark \ref{meas}.
To prove the last assertion of point 2, just notice that the 
following PCAF of $u$:
\[
t\mapsto \int_0^t |u(s,\theta)-a| \, l^a(ds,\theta),
\]
has Revuz measure:
\[
\sqrt{\frac 2{\pi\theta^3(1-\theta)^3}} \ a^2
\, \exp\left(-\frac{a^2}{2\theta(1-\theta)}\right) 
\, \cdot 
|x(\theta)-a| \, \nu(\, dx \,| \, x(\theta)=a) \, \equiv \, 0.
\]
To prove point 3 it is enough to notice that the PCAF
of $u$ in the left hand-side of (\ref{otfr})
has 1-potential $R(1)F_\theta$, where $F_\theta(y):=F(y(\theta))$,
$y\in K_0\cap C_0$, and the PCAF in the right hand side has
1-potential:
\begin{eqnarray*}
& &
\int_0^\infty r^2 \, F(r) \, dr
\int_{{\mathbb S}^{2}} \Pi \left( U_3^{\theta,r\cdot
n}\right)  \, {\cal H}^{2}(dn)
\, = \, \Pi\left(\int_{{\mathbb R}^3} F(|\alpha|) 
\, U_3^{\theta,\alpha} \, d\alpha \right)
\\ \\ & & = \, \Pi \left(R_3(1)[F_\theta\circ\Phi_3]\right) \, = \,
R(1)F_\theta.
\end{eqnarray*}
Since $R(1)F_\theta$ is bounded, then, arguing like in Theorem 5.1.6
of \cite{fot}, the two processes in (\ref{otfr})
coincide as PCAF's in the strict sense. \quad $\square$

\section{The reflecting measure $\eta$}

Recall that $\eta$ is the reflecting measure on
${\cal O}=[0,\infty)\times(0,1)$ which appears 
in equation (\ref{1}). The main result of this section
is the following:
\begin{theorem}\label{renor}
Let $\theta\in(0,1)$, $a\geq 0$.
\begin{itemize}
\item[1.] For all $\theta\in(0,1)$, there exists a PCAF in the
strict sense $(l(t,\theta))_{t\geq 0}$ of $u$, such
that $(l(\cdot,\theta))_{\theta\in(0,1)}$
is continuous in the sense of Definition \ref{def1}
and jointly measurable, and such that:
\[
l(t,\theta)
\, = \, \, \lim_{\epsilon\downarrow 0}
\frac 3{\,\epsilon^3} 
\int_0^t 1_{[0,\epsilon]}(u(s,\theta)) \, ds,
\]
in the sense of Definition \ref{def1}.
\item[2.] The PCAF $({l}(t,\theta))_{t\geq 0}$ has Revuz measure:
\[
\sqrt{\frac 2{\pi\theta^3(1-\theta)^3}} \, 
\, \nu(\, dx \,| \, x(\theta)=0),
\]
and increases only
on $\{t:u(t,\theta)=0\}$.
\item[3.] We have:
\[
{l}(t,\theta) \, = \, 
\lim_{a\downarrow 0} \frac 1{\, a^2} \, l^a(t,\theta)
\]
in the sense of Definition \ref{def1}.
\item[4.] For all $t\geq 0$ and $x\in K_0$, $\eta([0,t],d\theta)$ is absolutely
continuous w.r.t. the Lebesgue measure $d\theta$ and:
\begin{equation}\label{abscon.}
\eta([0,t],d\theta) \, = \, \frac 14 \, {l}(t,\theta) \, d\theta.
\end{equation}
\item[5.] For all $a\in (0,1)$:
\[
\lim_{\epsilon\downarrow 0}\, \sqrt\epsilon \,
\int_0^a \left(1 \wedge \frac{\theta}\epsilon \right)
\, \eta([0,t], d\theta) \, = \,
\sqrt{\frac 2\pi} \, t,
\]
\[
\lim_{\epsilon\downarrow 0}\, \sqrt\epsilon \,
\int_a^1 \left(1 \wedge \frac{1-\theta}\epsilon \right)
\, \eta([0,t], d\theta) \, = \,
\sqrt{\frac 2\pi} \, t
\]
in the sense of Definition \ref{def1}.
\end{itemize}
\end{theorem}
{\bf Proof.} For all $\epsilon>0$ we set:
\[
g^\epsilon(y):=
\frac 3{\,\epsilon^3} \,
1_{[0,\epsilon]}(y(\theta)),
\quad y\in K_0\cap C_0.
\]
By Lemma \ref{projector}, we have that:
\begin{eqnarray*}
R(1)g^\epsilon & = &
\Pi \left( R_3(1)\left[g^\epsilon\circ\Phi_3\right]
\right)
\, = \, \frac 3{\,\epsilon^3}
\int_{(|\alpha|\leq \epsilon)}
\Pi \left( U^{\theta,\alpha}_3 \right)
\, d\alpha
\\ \\ & = & 
\frac 3{\,\epsilon^3} \int_0^{\epsilon} r^2 \, dr
\int_{{\mathbb S}^{2}} \Pi \left( U_3^{\theta,r\cdot
n}\right)  \, {\cal H}^{2}(dn).
\end{eqnarray*}
By Lemma \ref{projector}, $R(1)g^\epsilon$ converges in
$W^{1,2}(\nu)$ as $\epsilon\to 0$ to:
\[
U^{\theta} \, := \, 4\pi \,
\Pi\left( U_3^{\theta,0}\right), \qquad a\geq 0.
\]
For all $\epsilon>0$ and
$\varphi\in W^{1,2}(\nu)\cap C_b(K_0)$ we have:
\[
{\cal E}_1(R(1)g^\epsilon,\varphi) \, = \, 
\int_{K_0} g^\epsilon \, \varphi \, d\nu \, = \, 
\frac 3{\,\epsilon^3} \,
{\mathbb E}\left[ \varphi(e) \, 1_{[0,
\epsilon]}(e(\theta))\right],
\]
and letting $\epsilon\to 0$ we get:
\begin{eqnarray}\label{smooth2}
\nonumber & &
{\cal E}_1(U^{\theta},\varphi) \, = \, \lim_{\epsilon\to 0}
\frac 3{\,\epsilon^3} \, {\mathbb E}\left[ \varphi(e) \, 1_{[0,
\epsilon]}(e(\theta))\right]
\\ \nonumber \\ & & = \,
\sqrt{\frac 2{\pi\theta^3(1-\theta)^3}}
\, {\mathbb E}\left[ \varphi(e) \, | \, e(\theta)=0\right].
\end{eqnarray}
Since $\Pi$ is Markovian, by (\ref{unibo}) $U^{\theta}$ is
bounded, and by (\ref{smooth2}) $U^{\theta}$ is the 
1-potential of a non-negative finite measure. By Theorem 5.1.6
in \cite{fot}, there exists a PCAF $(l(t,\theta))_{t\geq
0}$ in the strict sense of $u$, with 1-potential equal to
$U^{\theta}$ and with Revuz-measure given by
(\ref{smooth2}). Since $R(1)g^\epsilon$ is the 1-potential
of the following PCAF of $u$:
\[
t\mapsto \frac 3{\,\epsilon^3} 
\int_0^t 1_{[0,\epsilon]}(u(s,\theta)) \, ds, \quad t\geq 0,
\]
then, points 1 and 2 of Theorem \ref{renor} are proved
by (\ref{unifa}), Lemma \ref{lemdef} and Remark \ref{meas}.
To prove the last assertion of point 2, just notice that the 
following PCAF of $u$:
\[
t\mapsto \int_0^t u(s,\theta) \, l(ds,\theta),
\]
has Revuz measure:
\[
\sqrt{\frac 2{\pi\theta^3(1-\theta)^3}} \, 
\, \cdot x(\theta) \, \nu(\, dx \,| \, x(\theta)=0) \, \equiv \, 0.
\]
From the Proof of Theorem \ref{occuref}, we know that the
1-potential of $l^a(\cdot,\theta)$ is:
\[
U^{\theta,a} \, = \, a^2
\int_{{\mathbb S}^{2}}
\Pi\left( U_3^{\theta,a\cdot n}\right) \,
{\cal H}^{2}(dn), \qquad a\geq 0.
\]
Then $U^{\theta,a}/a^2$ converges as $a\to 0$
to $U^\theta$ in $W^{1,2}(\nu)$.
Since $U^{\theta,a}/a^2$ is the 1-potential
of $l^a(\cdot,\theta)/a^2$, by Lemma \ref{lemdef}
point 3 of Theorem \ref{renor} is proved.
Let now $I\subset\!\subset (0,1)$ be Borel. 
Notice that the following PCAF in the strict sense of $u$:
\begin{equation}\label{pcaf1}
t \, \mapsto \, \frac 14 \, \int_I l(t,\theta) \, d\theta
\end{equation}
has Revuz measure:
\begin{equation}\label{revin}
\frac 12 \, 
\int_I \frac 1{\sqrt{2\pi\theta^3(1-\theta)^3}} \, 
\, \nu(\, dx \,| \, x(\theta)=0) \, d\theta,
\end{equation}
and 1-potential equal to:
\[
\frac 14 \, U^I \, := \, \frac 14 \, \int_I U^\theta \, d\theta.
\]
On the other hand, it was proved in Theorem 7 of \cite{za00} 
that the PCAF in the strict sense of $u$:
\begin{equation}\label{pcaf2}
t \, \mapsto \, \eta([0,t]\times I) 
\end{equation}
has Revuz measure equal to (\ref{revin}). Therefore, 
by Theorem 5.1.6 in \cite{fot}, the two PCAFs of $u$ 
in (\ref{pcaf1}) and (\ref{pcaf2}) coincide, and 
since $U^I$ is a bounded 1-potential then they coincide
as PCAFs in the strict sense. Therefore point 4 is proved.
We prove now the last assertion.
For all $\epsilon\in(0,1/2)$
set $h_\epsilon:[0,1]\mapsto [0,1]$:
\[
h_\epsilon(\theta) \, := \, {\sqrt\epsilon}\left(
\left(1 \wedge \frac\theta\epsilon \right) 1_{[0,1/2]}(\theta)
\, + \, 4\theta(1-\theta) \, 1_{[1/2,1]}(\theta)\right).
\]
Then $h_\epsilon$ is concave and continuous on $[0,1]$, with:
\[
h_\epsilon''(d\theta) \, = \, - \, \frac 1{\sqrt\epsilon}
\, \delta_\epsilon(d\theta) \, - \, {\sqrt\epsilon} \, 8 \, 
1_{[1/2,1]}(\theta) \, d\theta,
\]
where $\delta_\epsilon$ is the Dirac mass at $\epsilon$. Moreover
$h_\epsilon(0)=h_\epsilon(1)=0$
and $h_\epsilon \to 0$ uniformly on $[0,1]$ as $\epsilon
\to 0$. By (\ref{weeq}) we have then:
\begin{equation}\label{limfi}
\lim_{\epsilon\to 0} \left(
\frac 1{2\sqrt\epsilon} \, \int_0^t u(s,\epsilon) \, ds \,
- \, \sqrt{\epsilon}\, \int_0^{1/2}
\left(1 \wedge \frac\theta\epsilon \right) \,
\eta([0,t],d\theta) \right) \, = \, 0.
\end{equation}
Recall the definition of $\overline{\gamma}^\theta$ given
in point 2 of Proposition \ref{preli}. 
We set $\gamma^\epsilon:K_0\cap C_0 \mapsto {\mathbb R}$,
$\gamma^\epsilon(x):=x(\epsilon)/{\sqrt\epsilon}$. Then, by
Lemma \ref{projector} we have that 
$R(1) \gamma^\epsilon \, = \, \Pi \left(R_3(1)
\overline{\gamma}^\epsilon\right)$.
By point 2 of Proposition \ref{preli} and by Lemma \ref{projector},
we obtain that $R(1) \gamma^\epsilon$ converges to 
$\sqrt{8/\pi}$ in $W^{1,2}(\nu)$. Therefore, by Lemma 
\ref{lemdef}:
\[
\lim_{\epsilon\to 0} 
\frac 1{2\sqrt\epsilon} \, \int_0^t u(s,\epsilon) \, ds \,
= \, \sqrt{\frac 2\pi} \, t,
\]
in the sense of Definition \ref{def1}, and by (\ref{limfi})
point 5 is proved.
\quad $\square$

\begin{corollary}\label{cor1}
For all $x\in K_0\cap C_0$, a.s. the set:
\[
S:=\{s>0: \, \exists \, \theta\in(0,1), \, u(s,\theta)=0\}
\]
is dense in ${\mathbb R}^+$ and has zero Lebesgue measure.
\end{corollary}
{\bf Proof.} By Point 5 in Theorem \ref{renor}, for all
$x\in K_0\cap C_0$, a.s. for all $t>0$ we have
$\eta([0,t]\times(0,1))=+\infty$, so that in particular
$\eta([0,t]\times(0,1))>0$. By (iv) in
Definition \ref{d1} the support of $\eta$ is contained in
the set $\{u=0\}$, so that for all $t>0$ there exists
$s\in(0,t)\cap S$. By the Markov property, for all
$q\in{\mathbb Q}$ and all $t>q$, there exists $s\in(q,t)
\cap S$, which implies the density of $S$ in
${\mathbb R}^+$. To prove that $S$ has zero Lebesgue measure,
recall that the law of $u(t,\cdot)$ is absolutely continuous
w.r.t. $\nu$ for all $t>0$, and $\nu(x:\exists \, \theta\in(0,1), \,
x(\theta)=0)=0$. Then, if ${\cal H}^1$ is the Lebesgue measure 
on ${\mathbb R}$:
\[
{\mathbb E}_x\left[{\cal H}^1(S)\right] = \int_0^\infty
{\mathbb E}_x\left[ 1_S(t) \right] \, dt =
\int_0^\infty {\mathbb P}(\exists \, \theta\in(0,1), \,
u(t,\theta)=0) \, dt = 0. \quad \square
\]

\vspace{.3cm}\noindent
Notice now that, by Points 2 and 4 of Theorem \ref{renor},
equation (\ref{1}) can be formally written in the following form:
\begin{equation}\label{11}
\left\{ \begin{array}{ll}
{\displaystyle
u(t,\theta)=x(\theta) + \frac 12 \int_0^t
\frac{\partial^2 u}{\partial \theta^2}(s,\theta) \, ds +
\frac{\partial W}{\partial\theta}(t,\theta) +
\frac 14 \, l(t,\theta) }
\\ \\
u(t,0)=u(t,1)=0
\\ \\
{\displaystyle
u\geq 0, \ l(dt,\theta)\geq 0, \
\int_0^\infty u(t,\theta)\, l(dt,\theta) =0, \quad \forall
\theta\in(0,1), }
\end{array} \right.
\end{equation}
where, as usual, the first line is rigorously defined
after taking the scalar product in $H$ between each term 
and any $h\in D(A)$.
Formula (\ref{11}) allows to interpret
$(u(\cdot,\theta),l(\cdot,\theta))_{\theta\in(0,1)}$ as
solution of a system of 1-dimensional Skorohod problems,
parametrized by $\theta\in(0,1)$. This fact is reminiscent of
the result of Funaki and Olla who proved in \cite{fuol}
that the stationary solution of a certain system of 1-dimensional 
Skorohod problems converges under a suitable rescaling
to the stationary solution of (\ref{1}).

\vspace{.3cm}\noindent
Finally, we show that $u$ satisfies a closed formula and that
equation (\ref{1}) is related to a fully non-linear equation.
Let $(w(t,\theta))_{t\geq 0, \theta\in[0,1]}$ be the
Stochastic Convolution:
\[
w(t,\theta):=\int_0^t\int_0^1 g_{t-s}(\theta,\theta') \, 
W(ds,d\theta'),
\]
solution of:
\begin{equation}\label{oud1}
\left\{ \begin{array}{ll}
{\displaystyle
w(t,\theta)= \frac 12 \int_0^t
\frac{\partial^2 w}{\partial \theta^2}(s,\theta) \, ds +
\frac{\partial W}{\partial\theta}(t,\theta) }
\\ \\
w(t,0)=w(t,1)=0
\end{array} \right.
\end{equation}
Subtracting the first line of (\ref{11}) and the first line 
of (\ref{oud1}), we obtain that:
\[
(t,\theta) \, \mapsto \,
\frac 12 \frac{\partial^2}{\partial \theta^2}
\int_0^t \left(u(s,\theta)-w(s,\theta)\right) \, ds
\]
is in $L^1_{loc}((0,1);$ $C([0,T]))$ for all
$T>0$, i.e. admits
a measurable version which is continuous in $t$
for all $\theta\in(0,1)$ and such that 
the sup-norm in $t\in[0,T]$ is locally integrabile
in $\theta$. Then, we can write:
\begin{equation}\label{111}
\left\{ \begin{array}{ll}
{\displaystyle
u(t,\theta)=x(\theta) + w(t,\theta) +
\frac 12 \frac{\partial^2}{\partial \theta^2}
\int_0^t \left(u-w\right)(s,\theta) \, ds + 
\frac 14 \, l(t,\theta)
}
\\ \\
u(t,0)=u(t,1)=0
\\ \\
{\displaystyle
u\geq 0, \ l(dt,\theta)\geq 0, \
\int_0^\infty u(t,\theta)\, l(dt,\theta) =0, \quad \forall
\theta\in(0,1), }
\end{array} \right.
\end{equation}
where every term is now well-defined and continuous in $t$,
and we can apply Skorohod's Lemma (see Lemma VI.2.1 in \cite{reyo})
for fixed $\theta\in(0,1)$, obtaining:
\begin{equation}\label{sk}
\frac 14 \, l(t,\theta) \, = \,
\sup_{s\leq t} \left[ -\left( x(\theta) + 
w(s,\theta) \, + \, \frac 12 
\frac{\partial^2}{\partial \theta^2}  
\int_0^s \left(u-w\right)(r,\theta) \, dr \right)\right] \vee 0,
\end{equation}
for all $t\geq 0$, $\theta\in(0,1)$. Therefore we have the
following:
\begin{corollary}\label{cor2}
For all $x\in K_0\cap C_0$, a.s. $u$ satisfies the
closed formula:
\begin{eqnarray}\label{eqfu}
u(t,\theta) & = & x(\theta) + w(t,\theta) \, + \, \frac 12 
\frac{\partial^2}{\partial \theta^2} 
\int_0^t \left(u(s,\theta)-w(s,\theta)\right) \, ds 
\\ \nonumber \\ \nonumber
& + & \sup_{s\leq t} \left[ -\left( x(\theta) + 
w(s,\theta) \, + \, \frac 12 
\frac{\partial^2}{\partial \theta^2}  
\int_0^s \left(u(r,\theta)-w(r,\theta)
\right) \, dr \right)\right] \vee 0,
\end{eqnarray}
for all $t\geq 0$, $\theta\in(0,1)$.
In particular $v$, defined by:
\[
v(t,\theta) \, := \, \int_0^t \left( u(s,\theta) -
w(s,\theta) \right) \, ds,
\]
is solution of the following fully non-linear equation:
\begin{equation}\label{fully}
\left\{ \begin{array}{ll}
{\displaystyle
\frac{\partial v}{\partial t}=
\frac 12 \frac{\partial^2 v}{\partial \theta^2} + x(\theta)
+ \sup_{s\leq t} \left[ -\left( x(\theta) + 
w(s,\theta) \, + \,
\frac 12 \frac{\partial^2 v}{\partial \theta^2}(s,\theta) 
\right) \right] \vee 0 }
\\ \\
v(0,\cdot)=0, \quad v(t,0)=v(t,1)=0,
\end{array} \right.
\end{equation}
unique in $\ {\cal V}:=\{v':\overline{\cal O}\mapsto
{\mathbb R}$ continuous: $\partial v'/\partial t$
continuous, $\partial^2 v'/\partial \theta^2
\in L^1_{loc}((0,1);$ $C([0,T]))$ for all $T>0\}$. 
\end{corollary}
The uniqueness of solutions of equation (\ref{fully}) in ${\cal V}$
is a consequence of the pathwise
uniqueness of solutions of equation (\ref{1}),
proved in \cite{nupa}: indeed, if $v'\in{\cal V}$ is a solution of 
(\ref{fully}), then setting
\[
u'(t,\theta) \, := \, \frac{\partial v'}{\partial t}(t,\theta)
\, + \, w(t,\theta), 
\]
\[
\eta'(dt,d\theta) \, := \, d_t \left\{
\sup_{s\leq t} \left[ -\left( x(\theta) + 
w(s,\theta) \, + \,
\frac 12 \frac{\partial^2 v'}{\partial \theta^2}(s,\theta) 
\right) \right] \vee 0 \right\} \, d\theta
\]
and repeating the above arguments backwards, we obtain that
$(u',\eta')$ is a weak solution of (\ref{1}), so that $u'=u$
and therefore $v'=v$. 

Notice that, by point 5 in Theorem 
\ref{renor}, by (\ref{abscon.})-(\ref{sk})-(\ref{fully})
and by the continuity of $\partial v/\partial t$ on
$\overline{\cal O}$, then,
for all $t>0$, $\partial^2 v/\partial \theta^2(t,\cdot)$ is 
not in $L^1(0,1)$, so that by the uniqueness 
a $C^{1,2}(\overline{\cal O})$
solution of (\ref{fully}) does not exist.

\end{document}